\documentclass[review]{elsarticle}

\usepackage{lineno,hyperref}
\modulolinenumbers[5]
\usepackage{amsmath}
\usepackage{amsfonts}
\usepackage{amssymb}  
\usepackage[amsmath, thmmarks]{ntheorem}

\makeatletter
\newtheoremstyle{nummermitklammern}%
 {\item[\rlap{\vbox{\hbox{\hskip\labelsep \theorem@headerfont
  (##2)\ ##1\theorem@separator}\hbox{\strut}}}]}%
 {\item[\rlap{\vbox{\hbox{\hskip\labelsep \theorem@headerfont
  (##2)\ ##1 (##3)\theorem@separator}\hbox{\strut}}}]}%
\makeatother
\theoremstyle{nummermitklammern}
\theorembodyfont{\rmfamily}
\theoremsymbol{\ensuremath{\diamond}}
\newtheorem{defsatzusw}{}[section]
\newtheorem{definition}[defsatzusw]{Definition}

\newtheorem{theorem}[defsatzusw]{Theorem}
\newtheorem{lemma}[defsatzusw]{Lemma}

\newtheorem{proposition}[defsatzusw]{Proposition}
\newtheorem{remark}[defsatzusw]{Remark}
\newtheorem{example}[defsatzusw]{Example}
\newtheorem{notation}[defsatzusw]{Notation}
\theoremstyle{nonumberbreak}
\theoremsymbol{\ensuremath{\square}}
\newtheorem{proof}{Proof}

\newcommand{\Q}{\mathbb{Q}}

\newcommand{\C}{\mathbb{C}}

\journal{Bulletin des Sciences Math\'{e}matiques}

\usepackage{numcompress}

\begin{document}

\begin{frontmatter}

\title{Invariant ideals of local analytic and formal vector fields}
\tnotetext[mytitlenote]{\text{This work was supported by DFG Graduiertenkolleg 1632 "$\text{Experimental and constructive algebra}$"}}

\author{Niclas Kruff}
\address{Lehrstuhl A f\"{u}r Mathematik\\RWTH Aachen University\\ Templergraben 55, 52056 Aachen, Germany}

\ead{niclas.kruff@matha.rwth-aachen.de}



\begin{abstract}
The objective of this paper is to analyse analytic invariant sets of analytic ordinary differential equations (ODEs). For this purpose we introduce semi-invariants and invariant ideals as well as the notion of vector fields in Poincar\'{e}-Dulac normal form (PDNF). We prove that all invariant ideals of a vector field in PDNF are already invariant for its semi-simple linear part. Additionally, this paper provides a natural characterization of invariant ideals via semi-invariants.   
\end{abstract}


\end{frontmatter}


\section{Introduction and preliminaries}

\noindent In this paper we focus on analytic invariant sets of an ODE system
\begin{equation}\label{Equation}
 \dot x=f(x),\ x\in \mathbb{K}^n,\ \mathbb{K}\in\{\mathbb{R},\mathbb{C}\},\ n\in \mathbb{N},
\end{equation}
where $f$ is analytic. 
Besides being of interest in its own right, the discussion of local analytic invariant sets is also relevant for the discussion of invariant algebraic curves of planar polynomial vector fields considering \textit{stationary points at infinity}. For details see \cite{schlo1}, \cite{schlo2} and \cite{walcher1}.\\
\noindent If $f$ is analytic near a non-stationary point $x_0$, i.e.\ $f(x_0)\neq 0$, then we can apply the straightening theorem which clarifies the structure of analytic invariant sets near $x_0$. Consequently, we examine invariant sets near a stationary point of $f$, which can be assumed to be $0$. In order to find invariant sets near stationary points one can investigate \textit{semi-invariants} which are analytic functions that remain invariant under the action of the \textit{Lie derivative} $L_f$ ($L_f$-invariant) and which induce invariant sets of codimension one. A natural way to generalize the notion of semi-invariants is to deal with \textit{invariant ideals} over the ring of convergent power series about zero $\mathbb{K}\{x\}$. In the case $\mathbb{K}=\mathbb{C}$ invariant radical ideals have a one to one correspondence to anayltic invariant sets.\\
\noindent A common tool to study the local dynamcis near the stationary point $0$ is the Poincar\'{e}-Dulac normal form of $f$ (PDNF), characterized by the condition\ \[[f,B_s](x)=B_s\cdot f(x)-Df(x)\cdot B_s x=0\] where $B_s$ is the semi-simple part of $Df(0)$. Given an analytic invariant set we associate with it its invariant radical ideal. By transforming a vector field into PDNF the transformed invariant ideal remains invariant for the normal form. However, studying vector fields in PDNF requires to work in the ring of formal power series and therefore there may not exist a one to one correspondence between the transformed invariant ideal and analytic invariant sets. But we still have the following picture:
\begin{align*}
 V\subseteq \mathbb{K}^n\ \text{invariant set for (\ref{Equation})}\ \longleftrightarrow&\ I\subseteq \mathbb{K}\{x\}\ L_f\text{-invariant ideal}\\
   &\ \ \ \ \ \ \ \ \ \ \Big\downarrow f\longrightarrow f^*\ \text{PDNF}\\
  & I^*\subseteq \mathbb{K}[[x]]\ L_{f^*}\text{-invariant ideal}
\end{align*}

\noindent Vector fields in PDNF admit a one-parameter group of symmetries given by $\exp(tB_s)$.
There are some further noteworthy properties of a vector field $f$ in PDNF:
\begin{enumerate}[$\bullet$]
 \item Commuting vector fields: $[f,g]=0\Longrightarrow [B_s,g]=0$ which is a consequence of Lemma 1.2 in \cite{walcher2} (see also \cite{bruno}).
 \item First integrals: $L_{f}(\varphi)=0\Longrightarrow L_{B_s}(\varphi)=0$ (\cite{walcher2} Proposition 1.8, \cite{zhang} Lemma 2.3).
 \item Semi-invariants: If $L_{f}(\varphi)=\lambda\cdot\varphi$, $\lambda\in \mathbb{K}[[x]]$, holds then there exists a unit $\mu$ and a constant $c$ such that $L_{B_s}(\mu\varphi)=c\cdot \mu\varphi$ (\cite{walcher1} Lemma 2.2.). 
\end{enumerate}
Especially, by identifying $\varphi$ with the ideal $\langle \varphi\rangle$ (unique up to multiplication with units), the last condition can be rewritten in the form 
\begin{equation}\label{propPDNF}
 L_{f}\left(\langle \varphi\rangle\right)\subseteq \langle \varphi\rangle \Longrightarrow L_{B_s}\left(\langle \varphi\rangle\right)\subseteq \langle \varphi\rangle.
\end{equation}
Thus every principal ideal which is invariant under $L_f$ is already invariant under $L_{B_s}$. The main result of this paper is a generalization of (\ref{propPDNF}) to arbitrary invariant ideals over the ring of formal power series. More precisely, given an $L_f$-invariant ideal over the ring of formal power series, we will prove that this ideal is $L_{B_s}$-invariant.\\
In summary, if we start with an analytic invariant set of (\ref{Equation}) we have a look at its corresponding $L_f$-invariant ideal $I$. After transforming $f$ into a vector field $f^*$ in PDNF we obtain a transformed invariant ideal $I^*$ of the normal form. However, in general we lose analyticity of the transformation and of $f^*$. As a consequence we work in the ring of formal power series and we get that $I^*$ is invariant for the semi-simple part of $Df^*$ in the stationary point. Our main results can be summarized in the following picture.
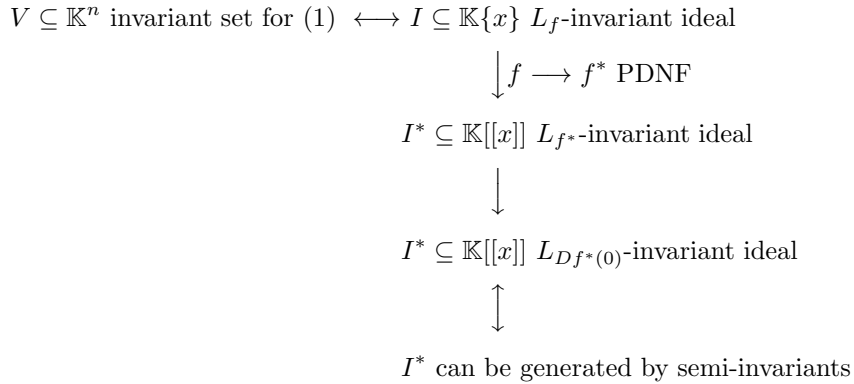
\begin{align*}
 V\subseteq \mathbb{K}^n\ \text{invariant set for (\ref{Equation})}\ \longleftrightarrow&\ I\subseteq \mathbb{K}\{x\}\ L_f\text{-invariant ideal}\\
  &\ \ \ \ \ \ \ \ \ \ \Big\downarrow f\longrightarrow f^*\ \text{PDNF}\\
  & I^*\subseteq \mathbb{K}[[x]]\ L_{f^*}\text{-invariant ideal}\\
  &\ \ \ \ \ \ \ \ \ \ \Big\downarrow\\
  &I^*\subseteq \mathbb{K}[[x]]\ L_{Df^*(0)}\text{-invariant ideal}\\
  &\ \ \ \ \ \ \ \ \ \ \Big\updownarrow\\
  &I^*\ \text{can be generated by semi-invariants}
\end{align*}
We close the paper with a class of examples.

\section{Notation and notions}

\noindent We first introduce some notations and recall some facts from the theory of vector fields in PDNF and from commutative algebra.

\begin{notation}
\noindent  We abbreviate $\mathbb{K}\{x_1,\dots,x_n\}=\mathbb{K}\{x\}$ for the ring of convergent power series about $0$ and $\mathbb{K}[[x_1,\dots,x_n]]=\mathbb{K}[[x]]$ for the ring of formal power series with $\mathbb{K}\in\{\mathbb{R},\mathbb{C}\}$. Moreover, we denote by $\text{Mon}(x)$ the set of all monomials.
 Given a vector field $f\in \mathbb{K}\{x\}^n$, $n\in \mathbb{N}$, with stationary point $0$, we have a look at the Taylor expansion about $0$
 \[
  f(x)=Bx+\sum\limits_{j\geq 2}^{}f^{(j)}(x)
 \]
 with $f^{(j)}$ homogeneous polynomial vector fields of degree $j$. The linear part $B\in \mathbb{K}^{n\times n}$ is given by $B=Df(0)$ which admits the \textit{Jordan-Chevalley decomposition} 
 \[
  B=B_s+B_n
 \]
into semi-simple part $B_s$ and nilpotent part $B_n$. The \textit{spectrum} $\text{Spec}(B)=\{\lambda_1\dots,\lambda_n\}\subseteq \mathbb{C}$ is given by all eigenvalues of $B$ counted with multiplicities.
 Given a function $\psi\in \mathbb{K}\{x\}$ its Lie derivative with respect to the vector field $f$ will be denoted by
 \[
  L_f(\psi)(x):=D(\psi)(x)\cdot f(x).
 \]
It describes the directional derivative of $\psi$ along solution trajectories of $\dot x=f(x)$. If $f$ is a formal vector field and $\psi$ a formal power series the formal Lie derivative of $\psi$ w.r.t. $f$ is still defined. We list up some properties of the Lie derivative in the appendix. Finally given two analytic or given two formal vector fields $f,g$ we define the Lie bracket of theses vector fields by
\[
 [f,g](x):=D(g)(x)\cdot f(x)-D(f)(x)\cdot g(x).
\]
One can easily verify the equation
\[
L_{[f,g]}=L_fL_g-L_gL_f. 
\]

\end{notation}

\begin{definition}
Let $f$ be a convergent or a formal vector field. Then $f$ is in PDNF if $[f,B_s](x)=0$. 
\end{definition}

\begin{remark}
By using Theorem 2.2 in \cite{bibikov}, for every vector field $f\in \mathbb{K}\{x\}^n$ there exists a near identity transformation into a vector field $\widehat{f}$ which is in PDNF. However, the transformation is not convergent in general and neither is $\widehat{f}$. For that reason we deal with vector fields over the ring of formal power series.
\end{remark}

\begin{definition}
Let $m\in\text{Mon}(x)$ with $m=x_1^{\alpha_1}\cdots x_n^{\alpha_n}$, $\alpha_j\in \mathbb{N}_0$, and $B_s=\text{diag}(\lambda_1,\dots,\lambda_n)$. We call 
\[
 w(m):=\sum\limits_{j=1}^{n}\alpha_j\lambda_j\in \mathbb{C}
\]
the weight of $m$. Since we can identify a monomial $m$ with its exponent we sometimes use the notation $w(\alpha)$ for $\alpha=\left(\alpha_1,\dots,\alpha_n\right)^{\text{tr}}$.\\
For a power series $\varphi$ we denote by $W(\varphi)$ the set of weights that appear in the monomial representation of $\varphi$.
\end{definition}

\noindent In particular we see that every monomial $m\in \text{Mon}(x)$ lies in an eigenspace of $L_{B_s}$ with $L_{B_s}(m)=w(m)\cdot m$.

\noindent We now come to the central objects in this paper.

\begin{definition}
\begin{enumerate}[(a)]
 \item Let $f\in \mathbb{K}\{x\}^n$ and $I\subseteq \mathbb{K}\{x\}$ be an ideal. We call the ideal $L_{f}$-invariant (or simply invariant) if $L_{f}(I)\subseteq I$ holds.
 \item In the special case of an invariant principal ideal which is generated by $\psi$, there exists $\lambda\in \mathbb{K}\{x\}$ such that $L_{f}(\psi)=\lambda\cdot \psi$. We then call $\psi$ a \textit{semi-invariant} of $f$.
 \hfill $\diamond$
\end{enumerate}
\end{definition}

\begin{remark}\label{nonewweights}
\begin{enumerate}[(a)]
\item The power series with a given weight form a subspace which is $L_{B_s}$-invariant.
\item  If $f$ is a vector field in PDNF and $m\in \text{Mon}(x)$ we observe the implication
 \[
  L_{B_s}(m)=w(m)m \Longrightarrow L_{B_s}(L_f(m))\overbrace{=}^{[f,B_s]=0}L_f(L_{B_s}(m))=w(m)L_f(m).
 \]
\noindent Hence, for $f$ in PDNF, the subspace of power series of a given weight is also $L_f$-invariant.
\end{enumerate}
\end{remark}

\noindent The study of invariant ideals is of particular interest because of their relation to invariant sets. For an ideal $I\subseteq \mathbb{K}\{x\}$ we define its vanishing set as
\[
 \mathcal{V}(I):=\{x\ |\ p(x)=0\ \forall\ p\in I\}.
\]
In the case $\mathbb{K}=\mathbb{C}$ invariant \textit{radical ideals} stand in one to one correspondence to invariant sets. (See \cite{chllpawa} Lemma 2.1. This is stated for polynomial ideals and in its proof one uses \textit{Hilbert's Nullstellensatz}. It can be generalized to the analytic case by using \textit{R\"{u}ckert's Nullstellensatz}.)

\noindent In the case $I\subseteq \mathbb{K}[[x]]$ we still talk about $L_{f}$-invariant ideals if $L_{f}(I)\subseteq I$ is satisfied. Consequently, an ideal $I\subseteq \mathbb{K}[[x]]$ is invariant if and only if $L_{f}(\varphi)\in I$ for all $\varphi\in I$.
In the following we investigate invariant ideals of formal vector fields in PDNF.

\section{Invariant ideals}

\noindent Throughout this section let $f\in K[[x]]^n$ be a formal vector field in PDNF, i.e.\ \[0=[f,B_s](x)=B_s\cdot f(x)-D(f)(x)\cdot B_sx,\] with diagonal semi-simple part $B_s=\text{diag}\left(\lambda_1,\dots,\lambda_n\right)$, $\lambda_i\in \mathbb{C}$. In the case $\mathbb{K}=\mathbb{R}$ we complexify if necessary.\\
\noindent The main goal of this paper is first to characterize all $L_{B_s}$-invariant ideals and secondly to generalize Lemma 2.2 in \cite{walcher1} (see (\ref{propPDNF})). The following two theorems are the main results of this paper.

\begin{theorem}\label{theo2}
Let $I\subseteq \mathbb{K}[[x]]$ be an $L_{B_s}$-invariant ideal. Then $I$ can be generated by semi-invariants of $B_s$, i.e.\ there exist $\varphi_1,\dots,\varphi_m\in \mathbb{K}[[x]]$ such that
\begin{enumerate}[(a)]
 \item $I=\langle \varphi_1,\dots,\varphi_m\rangle_{\mathbb{K}[[x]]}$.
 \item For all $1\leq j\leq m$ there exists $\nu_j\in \mathbb{K}[x]$ with the property $L_{B_s}(\varphi_j)=\nu_j\cdot \varphi_j$.
\end{enumerate}
\end{theorem}

\begin{theorem}\label{theo3}
 Let $I\subseteq K[[x]]$ be an $L_f$-invariant ideal. Then $I$ is $L_{B_s}$-invariant.
\end{theorem}

\noindent These theorems assert that semi-invariants of $B_s$ are the building blocks of invariant ideals, i.e.\ every invariant ideal can be generated by semi-invariants. We have divided the proofs of these theorems into a sequence of lemmas.

\begin{lemma}\label{Vandermonde}
 Let $I\subseteq \mathbb{K}[[x]]$ be an ideal. Assume that there exist \textit{polynomials} $\varphi_{i}\in \mathbb{K}[x]$, $1\leq i\leq l\in \mathbb{N}$, such that $I=\langle \varphi_{1},\dots,\varphi_{l}\rangle$.
 If $I$ is $L_{B_{s}}$-invariant then the ideal $I$ can be generated by semi-invariants of $B_s$.
\begin{proof}
Fix $1\leq i\leq l$ and let 
\begin{equation*}
 \varphi_{i}=\sum\limits_{k=1}^{r}c_{k}x^{\alpha_{k}},\ \alpha_{k} \text{ pairwise distinct},
\end{equation*}
where $r\in \mathbb{N}$ and $c_{k}\in \mathbb{C}\setminus \{0\}$. There exists a representation by the weights $w(\alpha_{k})$. Thus, we define the set
\begin{equation*}
W(\varphi_{i}):=\{w(\alpha_{k})\ |\ 1\leq k\leq r\}=:\{v_{1},\dots,v_{\mathfrak{q}}\}
\end{equation*}
where $\mathfrak{q}\leq r$ and $v_{k_{1}}\neq v_{k_{2}}$ whenever $k_{1}\neq k_{2}$ and we obtain a representation 
\begin{equation}\label{decomweights}
 \varphi_{i}=\sum\limits_{k=1}^{\mathfrak{q}}\left(\underbrace{\sum\limits_{1\leq m\leq r,\ w(\alpha_{m})=v_{k}}^{}c_{m}x^{\alpha_{m}}}_{\varphi_{i,v_{k}}}\right).
\end{equation}
By assumption, the ideal $I$ is $L_{B_{s}}$-invariant and this implies by induction
\begin{equation*}
 L_{B_{s}}^{(\ell)}(\varphi_{i})=\underbrace{L_{B_s}\circ\cdots\circ L_{B_s}}_{\ell\ \text{times}}(\varphi_i)\in I
\end{equation*}
for all $\ell\in \mathbb{N}_{0}$. Consequently, we have
\begin{equation*}
 L_{B_{s}}^{(\ell)}(\varphi_{i})=\sum\limits_{k=1}^{\mathfrak{q}}\left(\sum\limits_{1\leq m\leq r,\ w(\alpha_{m})=v_{k}}^{}v_{k}^{\ell}c_{m}x^{\alpha_{m}}\right)\in I
\end{equation*}
for all $0\leq \ell\leq \mathfrak{q}-1$. This induces the linear system of equations
\begin{equation*}
\underbrace{\begin{pmatrix} 1&1&\cdots &1 \\ v_{1}&v_{2}&\cdots &v_{\mathfrak{q}}\\ \vdots & \vdots & \cdots & \vdots\\ v_{1}^{\mathfrak{q}-1}&v_{2}^{\mathfrak{q}-1}&\cdots &v_{\mathfrak{q}}^{\mathfrak{q}-1}\end{pmatrix}}_{=:V}
\cdot \underbrace{\begin{pmatrix} \varphi_{i,v_{1}}\\ \varphi_{i,v_{2}}\\ \vdots \\ \varphi_{i,v_{\mathfrak{q}}} \end{pmatrix}}_{=:\Phi_{i}}
=\begin{pmatrix} L_{B_{s}}^{(0)}(\varphi_{i})\\ L_{B_{s}}^{(1)}(\varphi_{i})\\ \vdots \\ L_{B_{s}}^{(\mathfrak{q}-1)}(\varphi_{i}) \end{pmatrix}.
\end{equation*}
All components of the right hand side are elements of $I$. Moreover, the matrix $V\in \mathbb{K}^{\mathfrak{q}\times \mathfrak{q}}$ is a Vandermonde matrix, and hence it is invertible (see \cite{DaRe}, chapter 8). Therefore, all components of the vector $\Phi_{i}$ are elements of $I$, i.e.\ $\varphi_{i}$ has
a decomposition into semi-invariants with each semi-invariant an element of $I$. The index $i$ was arbitrary and the claim follows.
\end{proof}
\end{lemma}

\noindent We are now in a position to prove Theorem (\ref{theo2}). A crucial technique in order to prove to this theorem is to equip the ring $K[[x]]$ with the $x$-adic topology and to work in the quotient rings $K[[x]]/\langle x\rangle^i$, $i\in \mathbb{N}$ (see Remark (\ref{xadictopology})). 

\begin{proof}[of Theorem (\ref{theo2})]
Assume that $I\subseteq \mathbb{K}[[x]]$ is $L_{B_{s}}$-invariant and let $\varphi\in I$ arbitrary. The set
\begin{equation*}
W(\varphi):=\{w(x^{\alpha})\ |\ x^{\alpha}\text{ appears in the representation of }\varphi\}=:\{w_{k}\ |\ k\in Q\}, 
\end{equation*}
where $Q\subseteq\mathbb{N}$ and $w_k\neq w_{k'}$ for $k\neq k'$,
consists of all weights of monomials which appear in the representation of $\varphi$.
Analogously to Lemma (\ref{Vandermonde}) there exists a decomposition of the form
\begin{equation*}
 \varphi=\sum\limits_{w_{k}\in W(\varphi)}^{}p_{w_{k}}
\end{equation*}
where
\begin{equation*}
 p_{w_{k}}=\sum\limits_{\beta\in \mathbb{N}_{0}^{n},\ w(\beta)=w_{k}}^{}c_{\beta,w_{k}}x^{\beta},\ c_{\beta,w_{k}}\in \mathbb{C}\setminus \{0\}.
\end{equation*}
The cardinality of $W(\varphi)$ is infinite in general and the $p_{w_{k}}$ are formal power series. Due to the convergence in the $x$-adic topology (see \cite{GrePfi} Definition 6.1.6) we have to show that for all $\epsilon\in \mathbb{N}$ there exists a finite subset $E_{0}$ of $W(\varphi)$ such that
\[
 \left(\sum\limits_{w_{k}\in E}^{}p_{w_{k}}\right)-\varphi\in \langle x\rangle^{\epsilon}
\]
for all finite sets $E$ with the property
\[
 E_{0}\subseteq E\subseteq W(\varphi).
\]
However, this is obvious because if the monomial weights tend to infinity then the monomial exponents tend to infinity as well.
This yields that the sum
\begin{equation*}
 \sum\limits_{w_{k}\in W(\varphi)}^{}p_{w_{k}},
\end{equation*}
is a well-defined formal power series which equals $\varphi$.
It suffices to prove that $p_{w_{k}}\in I$ for all $w_{k}$. Consider the quotient ring 
\begin{equation*}
 R_{i}:=\mathbb{K}[[x]]/\langle x\rangle^{i},\ i\in \mathbb{N},
\end{equation*}
and $R_{0}:=\{0\}$.
By assumption, the ideal $I$ is $L_{B_{s}}$-invariant, i.e.
\begin{equation*}
 L_{B_{s}}(I)\subseteq I.
\end{equation*}
Consequently, one gets the relations
\begin{equation*}
 L_{B_{s}}((I+\langle x\rangle^{i})/\langle x\rangle^{i})=(L_{B_{s}}(I)+\langle x\rangle^{i})/\langle x\rangle^{i}\overbrace{\subseteq}^{(*)} (I+\langle x\rangle^{i})/\langle x\rangle^{i}.
\end{equation*}
The equivalence class of $\varphi$ with respect to $\langle x\rangle^{i}$ can be represented by a polynomial. By inclusion $(*)$ and Lemma (\ref{Vandermonde}), one has
\begin{equation*}
 [p_{w_{k}}]_{\langle x\rangle^{i}}\in (I+\langle x\rangle^{i})/\langle x\rangle^{i}
\end{equation*}
for all $i\in \mathbb{N}$. Hence, the residue class of $p_{w_k}$ is contained in the residue class of $I$ for all $i\in \mathbb{N}$.
This implies $p_{w_k}\in\overline{I}=I$ because ideals are closed sets under the $\langle x\rangle$-adic topology. In conclusion, the ideal $I$ can be generated by semi-invariants. By the Noetherian property of $\mathbb{K}[[x]]$ there are finitely many semi-invariants that generate $I$.
\end{proof}

\noindent In order to prove Theorem (\ref{theo3}) we need some additional considerations.

\begin{remark}
 Consider the decomposition
\begin{equation*}
 L_{f}=L_{B_{s}}+L_{g},
\end{equation*}
where $g=B_n+\sum\limits_{j\geq 2}^{}f^{(j)}$. The operator $L_{B_{s}}$ commutes with the operator $L_{g}$ due to the fact that $f$ is in PDNF, i.e.
\begin{equation*}
 [L_{B_{s}},L_{g}](\varphi):=L_{B_{s}}(L_{g}(\varphi))-L_{g}(L_{B_{s}}(\varphi))=0,\ \forall \varphi\in \mathbb{K}[[x]].
\end{equation*}
Hence, using the binomial theorem, one finds by induction
\begin{equation}\label{binomialformula}
 L_{f}^{(\ell)}(\varphi):=\left(\underbrace{L_{f}\circ \cdots \circ L_{f}}_{\ell\ \text{times}}\right)(\varphi)=\sum\limits_{j=0}^{\ell}\binom{\ell}{j}L_{B_{s}}^{(j)}(L_{g}^{(\ell-j)}(\varphi))
\end{equation}
for all $\ell\in \mathbb{N}_0$.\\
\noindent 
Consider a formal power series $\varphi\in \mathbb{K}[[x]]$. Motivated by the proof of Theorem (\ref{theo2}) we look at its residue class $[\varphi]_{R_i}$ in the ring $R_i:=\mathbb{K}[[x]]/\langle x^i\rangle$, $i\in \mathbb{N}$, which can be represented by a polynomial. Next, we investigate the decomposition w.r.t. to weights which is given by
\[
 [\varphi]_{R_i}=\sum\limits_{k=1}^{\mathfrak{q}}\left(\underbrace{\sum\limits_{1\leq m\leq r,\ w(\alpha_{m})=v_{k}}^{}c_{m}x^{\alpha_{m}}}_{\varphi_{v_{k}}}\right)\ \ \text{(compare to Equation (\ref{decomweights}))}.
\]
 Then we get:
\begin{align*}
\left[L_{f}^{(\ell)}([\varphi]_{R_i})\right]_{R_i}=\sum\limits_{k=1}^{\mathfrak{q}}\sum\limits_{j=0}^{\ell}L_{B_s}^{(j)}\left(\left[L_g^{(\ell-j)}(\varphi_{v_k})\right]_{R_i}\right).
\end{align*}
In the following we will concentrate on the term $\left[L_g^{(\ell-j)}(\varphi_{v_k})\right]_{R_i}$ and we will examine the structure of $L_{f}^{(\ell)}([\varphi]_{R_i})$.
\end{remark}

\noindent Due to Remark (\ref{nonewweights}), by applying $L_{f}$ to a formal power series $\varphi$ we know that all weights occurring in the monomial representation of $L_{f}(\varphi)$ are already contained in the set of weights of $\varphi$. So either some terms vanish or we get new terms with some weights that are already known, i.e.\ the set \[\{W(\left[L_f^{(\ell)}\left([\varphi]_{R_i}\right)\right]_{R_i})\ |\ \ell\in \mathbb{N}_0\}\] is finite and given by \[\{v_1,\dots,v_{\mathfrak{q}}\}.\]

\noindent This observation is crucial because we only have to consider the weights of the original element $\varphi$. The following lemma will be helpful.

\begin{lemma}\label{groworzero}
 Let $m\in \text{Mon}(x)$ be a monomial and $f$ a vector field in normal form. Furthermore let $f^{(j)}$ be a term in the representation of $f$ of degree $j\geq 2$. This implies
 \begin{equation*}
  \deg(m)< \deg(L_{f^{(j)}}(m)) \text{ or } L_{f^{(j)}}(m)=0.
 \end{equation*}
 \begin{proof}
  This is a direct consequence of
  \begin{equation*}
   \deg(L_{f^{(j)}}(m))=\deg(m)-j+1
  \end{equation*}
if $L_{f^{(j)}}$ does not annihilate $m$.
 \end{proof}
\end{lemma}

\noindent Furthermore, for all polynomials $\psi\in \mathbb{K}[x]$ there exists a $N=N(\psi)\in \mathbb{N}$ such that
\begin{equation*}
L_{B_{n}}^{(N)}(\psi)=0
\end{equation*}
because $B_{n}$ is a nilpotent matrix and therefore we get that $L_{B_{n}}|_{\mathcal{D}_k}$ is nilpotent for all $k\in \mathbb{N}$ (see Lemma (\ref{restrictionlie})).

\noindent Therefore, one gets the next proposition.

\begin{proposition}\label{nilpot}
 Let $\varphi\in \mathbb{K}[[x]]$ and $i\in \mathbb{N}$. There exists a natural number $\ell^{*}\in \mathbb{N}$ such that
 \begin{equation*}
  \left[L_{g}^{(\ell)}(\varphi)\right]_{R_i}=[0]_{R_{i}}\ \ ,
 \end{equation*}
 for all $\ell\geq \ell^{*}$, where $[\varphi]_{R_{i}}$ denotes the residue class of $\varphi$ in the ring $R_{i}$.
\end{proposition}

\noindent In conlusion there exists $\mathfrak{m}\in \mathbb{N}$ and a matrix $W\in \mathbb{C}^{\mathfrak{q}\cdot\mathfrak{m}\times \mathfrak{q}\cdot\mathfrak{m}}$ such that
\[
 L_{f}^{(\mathfrak{m})}(\varphi)=W\cdot\begin{bmatrix}L_{g}^{(0)}(\Phi)\\L_{g}^{(1)}(\Phi)\\ \vdots\\ L_{g}^{(\mathfrak{m}-1)}(\Phi)\end{bmatrix}
\]
with
\[
 \Phi=\begin{bmatrix}\varphi_{v_1}\\ \vdots\\ \varphi_{v_\mathfrak{q}}\end{bmatrix}
\]
and 
\[
 L_{g}(\Phi)=\begin{bmatrix}L_g(\varphi_{v_1})\\ \vdots\\ L_g(\varphi_{v_\mathfrak{q}})\end{bmatrix}.
\]

\noindent In general, by identity (\ref{binomialformula}), the matrix $W$ has the following structure:
\begin{equation*}
 W:=\begin{bmatrix}
    W_{1} & W_{2} & \cdots & W_{\mathfrak{m}}     
    \end{bmatrix}\in \C^{\mathfrak{q}\cdot\mathfrak{m}\times \mathfrak{q}\cdot\mathfrak{m}},
\end{equation*}
where $W_{i}\in \C^{\mathfrak{q}\cdot\mathfrak{m}\times \mathfrak{q}}$ and $\mathfrak{q},\mathfrak{m}\in \mathbb{N}$. The blocks $W_i$ are given by 
\begin{equation*}\begin{split}
 W_{1}&=\begin{bmatrix}
        1 & 1 & \cdots & 1\\
        w_{1} & w_{2} & \cdots & w_{\mathfrak{q}}\\
        \vdots & \vdots & \vdots & \vdots \\
        w_{1}^{\mathfrak{q}\mathfrak{m}-1} & w_{2}^{\mathfrak{q}\mathfrak{m}-1} & \cdots & w_{\mathfrak{q}}^{\mathfrak{q}\mathfrak{m}-1}
       \end{bmatrix}\\
 W_{2}&=\begin{bmatrix}
        0 & 0 & \cdots & 0\\
        1 & 1 & \cdots & 1\\
        \binom{2}{1}w_{1} & \binom{2}{1}w_{2} & \cdots & \binom{2}{1}w_{\mathfrak{q}}\\
        \vdots & \vdots & \vdots & \vdots \\
        \binom{\mathfrak{q}\mathfrak{m}-1}{1}w_{1}^{\mathfrak{q}\mathfrak{m}-2} & \binom{\mathfrak{q}\mathfrak{m}-1}{1}w_{2}^{\mathfrak{q}\mathfrak{m}-2} & \cdots & \binom{\mathfrak{q}\mathfrak{m}-1}{1}w_{\mathfrak{q}}^{\mathfrak{q}\mathfrak{m}-2}
       \end{bmatrix}\\ 
  W_{3}&=\begin{bmatrix}
        0 & 0 & \cdots & 0\\
        0 & 0 & \cdots & 0\\
        1 & 1 & \cdots & 1\\
        \binom{3}{2}w_{1} & \binom{3}{2}w_{2} & \cdots & \binom{3}{2}w_{\mathfrak{q}}\\
        \vdots & \vdots & \vdots & \vdots \\
        \binom{\mathfrak{q}\mathfrak{m}-1}{2}w_{1}^{\mathfrak{q}\mathfrak{m}-3} & \binom{\mathfrak{q}\mathfrak{m}-1}{2}w_{2}^{\mathfrak{q}\mathfrak{m}-3} & \cdots & \binom{\mathfrak{q}\mathfrak{m}-1}{2}w_{\mathfrak{q}}^{\mathfrak{q}\mathfrak{m}-3}
       \end{bmatrix}\\       
\end{split}\end{equation*}
and in general
\begin{equation*}
W_{p}=\begin{bmatrix}
       \binom{0}{p-1}w_{1}^{1-p} & \binom{0}{p-1}w_{2}^{1-p} & \cdots & \binom{0}{p-1}w_{\mathfrak{q}}^{1-p}\\
        \binom{1}{p-1}w_{1}^{2-p} & \binom{1}{p-1}w_{2}^{2-p} & \cdots & \binom{1}{p-1}w_{\mathfrak{q}}^{2-p}\\
        \binom{2}{p-1}w_{1}^{3-p} & \binom{2}{p-1}w_{2}^{3-p} & \cdots & \binom{2}{p-1}w_{\mathfrak{q}}^{3-p}\\
        \vdots & \vdots & \vdots & \vdots \\
        \binom{\mathfrak{q}\mathfrak{m}-1}{p-1}w_{1}^{\mathfrak{q}\mathfrak{m}-p} & \binom{\mathfrak{q}\mathfrak{m}-1}{p-1}w_{2}^{\mathfrak{q}\mathfrak{m}-p} & \cdots & \binom{\mathfrak{q}\mathfrak{m}-1}{p-1}w_{\mathfrak{q}}^{\mathfrak{q}\mathfrak{m}-p}\\
      \end{bmatrix}.
\end{equation*}

\begin{lemma}\label{genvandermonde}
The matrix $W$ is invertible.
\begin{proof}
 The matrix $W$ is a scaled transformation matrix of the Hermite interpolation\index{Hermite interpolation} problem with data points $w_{1},\dots,w_{l}$ and considering for each data point the first $\mathfrak{m}-1$ derivatives. The transformation matrix of the Hermite interpolation is also called \textit{confluent Vandermonde matrix} and is for example studied in \cite{dankal}.
 Since the Hermite interpolation problem has a unique solution (\cite{DaRe}, Satz 8.29) and therefore the matrix $W$ is invertible.
\end{proof}
\end{lemma}

\noindent Let us consider the following example to illustrate the argument.

\begin{example}
 Assume that $n=2$ and $f=\begin{bmatrix} x \\ 3y+\beta x^3\end{bmatrix}$ where $B_{s}=\text{diag}(1,3)$, $\beta\in \C$ and $f$ is in PDNF with $g=\begin{bmatrix}0\\ \beta x^3\end{bmatrix}$. Let $I\subseteq \C[[x,y]]$ be an $L_{f}$-invariant ideal and assume that $\psi:=x^3+y+y^2\in I$. By invariance, all expressions $L_{f}^{(\ell)}(\psi)$ are elements of $I$. Computing all those terms yields:
 \begin{equation*}\begin{split}
 &L_{f}^{(0)}(\psi)=\psi=3^{0}x^3+3^{0}y+6^{0}y^2\\
 &L_{f}^{(1)}(\psi)=\underbrace{3^{1}(x^3+y)+6^{1}y^2}_{=L_{B_{s}}(\psi)}+\underbrace{3^{0}\beta x^3+6^{0}2\beta x^3y}_{=L_{g}(\psi)}\\
 &L_{f}^{(2)}(\psi)=\underbrace{3^{2}(x^3+y)+6^{2}y^2}_{=L_{B_{s}}^{(2)}(\psi)}+2\cdot (\underbrace{3^{1}\beta x^3+6^{1}2\beta x^3y}_{=L_{B_{s}}(L_{g}(\psi))})+\underbrace{2\beta^{2} x^6}_{=L_{g}^{(2)}(\psi)}.
 \end{split}\end{equation*}
We observe that 
 \begin{equation*}
  L_{g}^{(k)}(\psi)=0
 \end{equation*}
holds for all $k\geq3$ and consequently this yields the following system of linear equations:
\begin{equation*}
 \underbrace{\begin{bmatrix}
  1 & 1 & 0 & 0 & 0 & 0\\
  3 & 6 & 1 & 1 & 0 & 0\\
  3^2 & 6^2 & 2\cdot 3 & 2\cdot 6 & 1 & 1\\
  3^3 & 6^3 & 3\cdot 3^2 & 3\cdot 6^2 & 3\cdot 3 & 3\cdot 6\\
  3^4 & 6^4 & 4\cdot 3^3 & 4\cdot 6^3 & 6\cdot 3^2 & 6\cdot 6^2\\
  3^5 & 6^5 & 5\cdot 3^4 & 5\cdot 6^4 & 10\cdot 3^3 & 10\cdot 6^3\\
 \end{bmatrix}}_{=:W}\cdot
 \begin{bmatrix}
  x^3+y\\
  y^2\\
  \beta\cdot x^3\\
  2\beta\cdot x^3y\\
  0\\
  2\beta^2\cdot x^6
 \end{bmatrix}\in I^{6}.
\end{equation*}
The matrix $W$ is invertible 
because 
\begin{equation*}
 \left|\det(W)\right|=(6-3)^9\neq 0
\end{equation*}
or by applying Lemma (\ref{genvandermonde}).
In conclusion, the polynomials $x^3+y$ and $y^2$ are elements of $I$.
\end{example}

\noindent Combining Proposition (\ref{nilpot}) and Lemma (\ref{groworzero}) we are finally in a position to prove Theorem (\ref{theo3}).

\begin{proof}[of Theorem (\ref{theo3})]
 Let $\varphi\in I$ and consider a decomposition $\varphi=\sum\limits_{\text{weights } w_k}^{}p_{w_k}$ where all $p_{w_k}$ are again semi-invariants of $B_{s}$. Fix $i\in \mathbb{N}$ and consider $[\varphi]_{R_{i}}$ which can be represented by a polynomial. 
 Using Lemma (\ref{nilpot}), there exists $\ell^{*}\in \mathbb{N}$ such that $L_{g}^{(\ell)}([\varphi]_{R_{i}})=0$ for all $\ell\geq \ell^{*}$. Again, one gets a (finite) linear system of equations and by Lemma (\ref{genvandermonde}) the associated matrix is invertible. 
 Thus, the linear system is uniquely solvable. Consequently, this yields
 \begin{equation*}
  [p_{w_{k}}]_{R_{j}}\in (I+\langle x\rangle^{j})/\langle x\rangle^{j},\ \forall j\in \mathbb{N}.
 \end{equation*}
Ideals are closed sets under the $x$-adic topology and therefore $p_{w_{k}}\in I$.
\end{proof}

\noindent As an application, we have a look at the following class of examples.

\begin{example}
Let $f\in \mathbb{K}[[x]]$ be a formal vector field in PDNF and $B_s=\text{diag}(\lambda_1,\dots,\lambda_n)$ such that $\lambda_1,\dots,\lambda_{n-1}$ are linearly independent over $\Q$ and 
\[
 \lambda_n=\sum\limits_{i=1}^{n-1}\alpha_i\lambda_i,\ \alpha_i\in \Q_{\leq 0}.
\]
This is the so called \textit{single resonance case}.
We want to determine all possible candidates of $L_{f}$-invariant prime ideals since all invariant ideals can be constructed by invariant prime ideals.\\
\noindent We first investigate all $L_{B_s}$-invariant ideals. By applying (\ref{theo2}) all invariant ideals can be genereated by semi-invariants and by using (\ref{propPDNF}) the only irreducible semi-invariants of $B_s$, up to multiplication with invertible power series, are given by 
\[
 x_1,\dots,x_n.
\]
To see this let $\varphi\in \mathbb{K}[[x]]$ be an irreducible semi-invariant of $B_s$. By applying part three of (\ref{propPDNF}), we may assume that there exists a constant $c\in \mathbb{K}$ such that
\[
 L_{B_s}(\varphi)=c\cdot \varphi.
\]
Hence, there exists a monomial $m=x_1^{\beta_1}\cdots x_n^{\beta_n}$ in the representation of $\varphi$ such that $w(m)=c$ and with $\beta_n\in \mathbb{N}_0$ minimal. Assuming that there exists an additional monomial $\widetilde{m}=x^{\gamma}\in \text{Mon}(x)$, with $w(\widetilde{m})=c$, implies, by using the linear independence of $\lambda_1,\dots,\lambda_{n-1}$, that 
\[
 \beta_j-\gamma_j=(\gamma_n-\beta_n)\alpha_j
\]
holds for all $1\leq j\leq n-1$. The numbers $\alpha_j$ are non-positive and $\beta_n$ is minimal and thus we obtain $\gamma_j\geq\beta_j$ for all $1\leq j\leq n$. In conclusion there exists a unit $u\in \mathbb{K}[[x]]^*$ such that 
\[
 \varphi=u\cdot x^{\beta}.
\]
Therefore, the irreducible semi-invariants, up to multiplication with units, are $x_1,\dots,x_n$.
Now, Theorem (\ref{theo3}) states that all $L_f$-invariant ideals are $L_{B_s}$-invariant ideals and these can be generated by semi-invariants. In conclusion all $L_f$-invariant prime ideals have the form
\[
 \langle x_{i_1},\dots,x_{i_r}\rangle, 1\leq r\leq n,
\]
with $i_1,\dots,i_r\in \{1,\dots,n\}$ and $i_1<i_2<\cdots<i_r$.
\end{example}

\section{Appendix}

\subsection{Commutative algebra}
\begin{remark}\label{xadictopology}
Let $R$ be a ring and $I\subseteq R$ an ideal. As shown in chapter 8 in \cite{mat}, the sets $f+I^{j}$, where $f\in R$ and $j\geq 1$, form a basis of open sets of the $I$-adic topology on $R$. If $S\subseteq R$ is an ideal the closure of $S$ is given by
\begin{equation*}
 \overline{S}:=\bigcap\limits_{j=1}^{\infty}(S+I^{j}).
\end{equation*}
 Now, let $J\subseteq \mathbb{K}[[x]]$ be an ideal. The closure of $J$, denoted by $\overline{J}$, is the set
 \begin{equation*}
  \overline{J}=\bigcap\limits_{i=1}^{\infty}(J+\langle x\rangle^{i}).
 \end{equation*}
However, we have the equality
\begin{equation*}
 \bigcap\limits_{i=1}^{\infty}\langle x\rangle^{i}=\{0\}
\end{equation*}
and this implies that each ideal $J\subseteq \mathbb{K}[[x]]$ is closed under the $\langle x\rangle$-adic topology.
\end{remark}

\subsection{Properties of the Lie derivative}

\begin{lemma}[see e.g.\ \cite{walchscheumay}, chapter 3]\label{restrictionlie}
Let $A\in \mathbb{K}^{n\times n}$ be a matrix with
\[
 \text{Spec}(A)=\{\lambda_{1},\dots,\lambda_{n}\}\subseteq \C
\]
and $\mathcal{D}_k$ be the space of homogeneous polynomials of degree $k$.
Moreover, we denote by \[L_{A}(\cdot):=L_{Ax}(\cdot)\] the Lie derivative with respect to the vector field $Ax$. Then we get:
\begin{enumerate}[i)]
\item $L_{A}(\mathcal{D}_k)\subseteq \mathcal{D}_{k}$, i.e.\ the restriction of $L_{A}$ to $\mathcal{D}_{k}$ is an endomorphism.
\item The eigenvalues of $(L_{A})|_{\mathcal{D}_{k}}$ are given by $\sum\limits_{i=1}^{n}m_{i}\lambda_{i}$ with
\[
(m_{1},\dots,m_{n})\in \mathbb{N}_{0}^{n}\ \text{such that}\  \sum\limits_{i=1}^{n}m_{i}=k.
\]
\item If $A=A_{s}+A_{n}$ is the decomposition into semi-simple\index{semi-simple} and nilpotent part\index{nilpotent} of $A$, then
 $(L_{A})|_{\mathcal{D}_{k}}=(L_{A_{s}})|_{\mathcal{D}_{k}}+(L_{A_{n}})|_{\mathcal{D}_{k}}$ is the decomposition into semi-simple and nilpotent part of $(L_{A})|_{\mathcal{D}_{k}}$.
\end{enumerate}
\end{lemma}

\end{document}